\documentclass[12pt]{amsart}
\usepackage{latexsym}
\usepackage{amsthm}
\usepackage{amsmath}
\usepackage{epsfig}
\usepackage{amsfonts}
\usepackage{amssymb}
\newtheorem{theorem}{Theorem}
\newtheorem{lemma}{Lemma}[theorem]
\newtheorem{defn}{Definition}
\newtheorem{prop}{Proposition}[theorem]

\newtheorem{example}{Example}

\def\R{\mathbb R}
\def\L{\mathcal{L}}

\begin{document}

\title{Total Curvature and Packing of Knots}

\author{Gregory Buck} \address{G. Buck, Department of Mathematics, St. Anselm
College,\newline Manchester, NH 03102.} 
\author{Jonathan Simon}\address{J. Simon, Department of Mathematics, University of Iowa, \newline Iowa City IA 52242.}
\thanks {Research supported by NSF Grants DMS\,0107747 (Buck) and DMS\,0107209 (Simon).  Email jsimon@math.uiowa.edu or gbuck@anselm.edu.  We thank J. McAtee and R. Weiler for helpful comments.}

\begin{abstract} We establish a new relationship between total
curvature of knots and crossing number. If $K$ is a smooth knot in
$\R^3$, $R$ the cross-section radius of a uniform  tube neighborhood
$K$,
$L$ the arclength of $K$, and $\kappa$ the total curvature of $K$, then

\begin{center}  crossing number of $K$ $< 4\,\frac{L}{R}\;
\kappa$\;.
\end{center}

The proof generalizes to show that for smooth knots in $\R^3$, the crossing
number, writhe, M\"obius Energy, Normal Energy, and Symmetric Energy are all bounded by the product of total curvature and rope-length.

One can construct knots in which the crossing numbers grow as fast as the
$(4/3)$ power of $\frac{L}{R}$. Our theorem says that such families must
have unbounded total curvature: If the total curvature is bounded, then the
rate of growth of crossings with ropelength can only be linear.

Our proof relies on fundamental lemmas about the total
curvature of curves that are packed in certain ways:  If a long smooth curve
$A$ with arclength $L$ is contained in a solid ball of radius $\rho$, then
the total curvature of $K$ is at least proportional to
$L/\rho$.  If $A$ connects concentric spheres of radii $a \geq 2$ and $b\geq
a+1$, by running from the inner sphere to the outer sphere and back again,
then the total curvature of $A$ is at least proportional to $1/\sqrt{a}$.
\end{abstract}

\maketitle

\markboth{G. Buck and J. Simon}{Total Curvature and Packing of Knots}

\section{Introduction} The total curvature of smooth closed  curve in
$\R^3$ must be at least $2\pi$; this is a theorem of Fenchel
\cite{Ch,Fen51}.  If the curve  actually is a nontrivial knot, then the
Fary-Milnor theorem  \cite{Ch,Fary, Fen51, MilnorTotalCurvature} says the
total curvature must be $> 4\pi$.  Are there properties of the knot that
could guarantee larger total curvature? Successive composition \cite{Fox50}
or or other kinds of satellite constructions (\cite{Schubert54} together
with \cite{MilnorTotalCurvature}) will work.  On the other hand, topological
complexity in the form of high crossing-number is not enough: it is well
known at least since \cite{Milnor53} that one can construct knots with
arbitrarily large minimum crossing-number  represented by curves with
uniformly bounded total curvature.  Here is one way to build them.
 
\begin{example}[Knots with bounded total curvature] Fix any odd integer
$n$.  Construct a smooth knot $K_n$, with minimum crossing number $n$, as
the union of four arcs $H,A,C_1, C_2$, where the total curvatures are
$\kappa(H) \to 0$ as $n \to \infty$,
$\kappa(A)=0$, and
$\kappa(C_1) \approx \kappa(C_2) \approx 2\pi$.  Let $H$ be the circular
helix in $\R^3$ parametrized as $[\cos(t), \sin(t), n^2t]$, $t=0\ldots
n\pi$.  {\rm The height coordinate $n^2t$ makes $\kappa(H)$ behave like
$1/n$ for large $n$.  Using any exponent larger than 1, i.e.
$n^{1+\epsilon}t$, still makes
$\kappa(H) \to 0$.}  Let $A$ be the central axis of the cylinder on which
$H$ runs.  Let $C_1$ and $C_2$ be curves that smoothly connect the top of
$H$ to the bottom of $A$ and vice-versa. {\rm For large $n$, the tangent
vectors at the beginning and end of $H$ are nearly vertical.  The arcs
$C_1$ and $C_2$ can be chosen to be almost planar-convex curves, with total
curvatures $\kappa(C_i) \approx 2\pi$.}  {\rm Similarly,} for any $(p,q)$,
torus knots or links of type ($p$ meridians,
$q$ longitudes), can have total curvature close to
$2\pi q$ if they are drawn on a standard torus that is long and thin enough.
\end{example}

In this paper, we show that examples of the preceeding kind are, in a sense,
the only kind possible.  In order to represent an infinite family of knot
types with uniformly bounded total curvature,  the knots must be ``long and
thin"; if we imagine them made of actual ``rope", then the ratio of length to
rope-thickness must grow without bound.

%%%%%%%%%%%%
%%%%%%%%%%%%
%If we form
%knots  by successive composition, e.g.
%$K\# K\ldots \# K$, then any
%curves representing these knot types will have the property that total
%curvature
%$\to \infty$ as the number of factors increases (\cite{Fox50} or
%\cite{Schubert54}, combined with
%\cite{MilnorTotalCurvature}).  By the same reasoning, if we form knots
%using some other kind of repeated satelite construction, then the bridge
%numbers, hence total curvature, grow without bound.|
%%%%%%%%%%%%
%%%%%%%%%%%%

\begin{defn} Suppose $K$ is a smooth knot in $\R^3$.  For $r>0$, consider
the disks of radius
$r$ normal to $K$, centered at points of $K$.  For $r$ sufficiently small,
these disks are pairwise disjoint and combine to form a tubular neighborhood
of $K$.  Let R(K), the {\em thickness radius} of $K$, denote the supremum of
such ``good" radii.  The {\em ropelength} of $K$, denoted
$E_L(K)$, is the ratio 
$$E_L(K)=  \frac{\text {\rm total arclength of } K}{R(K)}\;. $$
\end{defn}

The fundamental properties of thickness radius were developed in
\cite{LSDR}.  The idea of using the ratio of length--to--radius to measure
knot complexity was introduced in \cite{BO1}, and this ratio, denoted
$E_L(K)$, is connected to other knot ``energies" in \cite{BS2, BS3} and
\cite{RSBound}.  Variations on thickness are developed in 
\cite{Mo,devrthickness,DEJ2,devrthickness2, GonzalezMaddocks99,KS2,
RawdonIdealKnots98, Rawdon2000}.

\begin{defn} Let $K$ be a smooth knot.  From almost every direction, if we
project $K$ into a plane, the projection is regular, in particular there are
only finitely many crossings.  We can average this crossing-number over all
directions of projection (i.e. over the almost-all set of directions that
give regular projections).  This {\em average crossing number} is denoted
$\rm{acn}(K)$. 
  \rm {Certainly, the minimum crossing-number of the knot-type,}
$\rm{cr}[K]$, \rm {satisfies}
$\rm{cr}[K]
\leq
\rm{acn}(K)$.  We shall rely on the formulation of $\rm{acn}(K)$ developed
in \cite{FHW}.
\end{defn}

Our main result is the following: 
\begin{theorem}\label{MainTheorem} If $K$ is a smooth knot in $\R^3$, then
$$\rm{acn}(K) < 4\; E_L(K) \;\kappa(K)\;.$$
\end{theorem} 

The coefficient $c=4$ has been rounded up for simplicity.  What matters is that the crossing number is essentially bounded by ropelength times curvature.  In heuristic discussions, we may omit coefficients altogether.

If we are given some family of knots in which  total
curvature is uniformly bounded, while crossing number is growing, then the
ropelength must be growing at least as fast as the crossing numbers. 
Alternatively, if the crossing numbers are growing faster than ropelength,
then the total curvatures must be growing fast enough to make up the
difference.  We showed in  \cite{BS2,BS3} that $\rm{acn}(K)
\leq E_L(K)^{4/3}$, and there are examples \cite{B2, CKS} where the $4/3$
power is achieved.  In the particular examples of \cite{B2, CKS}, the knots
and links have evident growing total curvature; our theorem says that some
unbounded amount of total curvature must occur in any situation of
more-than-linear  growth of crossings with ropelength.

If we model a knot made of actual ``rope" as a smooth curve with a uniform
tube neighborhood, then the thickness (radius) $r$ of that rope is $\leq
R(K)$, so  $E_L(K) \leq \frac{L}{r}$.  Thus the theorem also holds with 
$\frac{L}{r}$ in place of $E_L$.

\section{Lemmas on total curvature}

The three lemmas in this section establish fundamental properties of smooth
space-curves, relating total curvature to packing, to oscillation relative
to a given point, and to the ``illumination'' of a given point.  We deal in
this section  with smooth space-curves, not assuming the curves are  simple
or closed; and we make no use of thickness.   Also we do not assume the
curves have finite length.

 To keep the arguments as simple as possible, we assume throughout the paper
that  ``smooth" means smooth of class
$C^2$.  The lemmas and theorem can be adapted for curves that are piecewise smooth.  For a smooth curve
$A$, we denote the total curvature of $A$ by $\kappa(A)$.

It is intuitively clear that if a long rope is packed in a small box, then
the rope must curve a lot.  This fundamental lemma is an important
ingredient in our analysis of the interplay between ropelength, crossing
number, and total curvature.  

A ball of radius $\rho$ contains, of course, a diameter of length $2\rho$. 
But once we postulate length $> 2\rho$, an arc in the ball must curve.  In
this version, we use $L\geq 3\rho$, but any constant larger than $2$ will
produce some guaranteed amount of total curvature.  Inequality (\ref{packinginequality1}) and the proof below are taken from \cite{chakerian64}, with a slight adjustment for non-closed curves.  If the curve is closed, then the number $2$ can be omitted from (\ref{packinginequality1}).

\begin{lemma}[Packing and curvature] \label{packing} Suppose $A$ is a smooth
connected curve  of length
$L$, contained in a round 3-ball of radius $\rho$.  Then
$\kappa(A)$ is approximately proportional to at least $L/\rho$.  More precisely, letting $\kappa$ denote $\kappa(A)$, we have the following,

\begin{equation}\label{packinginequality1}
L \leq \rho(\kappa +2)\;,\end{equation}
which gives 
\begin{equation} \label{packinginequality2} L \geq 3\rho \implies \kappa
\geq  1\;. \end{equation}  
\end{lemma}

\begin{proof}

Translate the ball and curve so the center of the ball is at the origin.  
Let $s \to x(s)$, $s \in [0,L]$, be a unit speed parametrization of $A$.  Since $|x'(s)|=1$, we can write 
$$L = \int_{s=0}^L x'(s) \cdot x'(s)\;.$$
Integrate by parts to get
$$L = x'(s)\cdot x(s) \mid_0^L - \int_{0}^L x(s)\cdot x''(s)\;.$$
Since $|x'(s)| = 1$, and $|x(s)| \leq \rho$, the first term is at most $2\rho$ and the second term is at most $\rho \kappa$.
\end{proof}

Another basic way that a long curve is forced to have a guaranteed amount of total
curvature is if its distance from some given point oscillates.  This is
captured in the next lemma.  

\begin{lemma}[Oscillation and curvature] \label{oscillation} Let $S_a, S_b$
be concentric spheres with radii $a<b$. Let $A$ be a smooth curve that
starts at a point of $S_a$, somewhere touches the sphere $S_b$, and ends at
a point of $S_a$.  Then the total curvature is at least approximately on the order of $1/\sqrt{a}$.

More precisely, 
 \begin{equation} \label{arcsininequality}
 \kappa(A) \geq \pi-2 \arcsin(a/b)\;. 
\end{equation}                If  $b \geq a+1$, then 
\begin{equation} \label{oscillationinequality}
\kappa(A) >\frac{ 2
\sqrt{2}}{\sqrt{a+1}}\;.
\end{equation} If $a \geq 2$, the bound (\ref{oscillationinequality}), along
with  simplifying the coefficient, gives  
$$
 \kappa(A) >\frac{2}{\sqrt{a}} \; . 
$$
\end{lemma}

{\bf Remark}.  One can think of the arc in this lemma as being just contained in the
spherical shell bounded by the two given spheres, or as reaching out past
$S_b$, just so it returns back to end on $S_a$.

\begin{proof} We rely on \cite{MilnorTotalCurvature} to reduce the proof of
(\ref{arcsininequality}) to analyzing a certain triangle, and then calculate
(\ref{oscillationinequality}).

Let
$p,q$ be the endpoints of
$A$, 
$z$  a point of
$A
\cap S_b$, 
 $\bar{pz}$ and $\bar{zq}$ the line segments from $p$ to $z$ and from $z$ to
$q$, and let $P$ be the two-edge polygon  $\bar{pz}  \cup \bar{zq}$.  Since
the polygon  $P$ is an inscribed polygon of $A$, we know from
\cite{MilnorTotalCurvature} that $\kappa(A) \geq \kappa(P)$, so it suffices
to establish the desired lower bound for $\kappa(P)$. If either of the edges
of $P$ is not tangent to the sphere $S_a$, we can pivot the edge at point
$z$ to move the edge to tangency in a way that opens the angle
$\widehat{pzq}$, so reducing $\kappa(P)$.  Thus, it suffices to prove the
lower bound for tangent two-edge polygons. 

The three points $p,z,q$ determine a plane; we wish that plane also would
include the center of the spheres.  If not, then (keeping the two edges
tangent to $S_a$, and allowing the points of tangency and the angle
$\widehat{pzq}$ to change), rotate the plane of
$p,z,q$ (with axis of rotation the line through $z$ parallel to
$\bar{pq}$) until it does contain the center of the spheres.  This
deformation also would increase the angle
$\widehat{pzq}$, and so decrease the total curvature of $P$. Thus, we are
reduced to the situation where $p,z,q$ and the center of the spheres are
coplanar, and $P$ consists of tangent lines to $S_a$ drawn symmetrically
from point $z$ on $S_b$.  We then have a right-triangle with
$$\sin (\frac{\widehat{pzq}}{2}) = \frac{a}{b}\;,$$ which gives
(\ref{arcsininequality}).

To derive (\ref{oscillationinequality}), first note that
$\pi-2 \text{arcsin }(a/b)$ increases as $(b-a)$ gets larger. So if we show 
$\pi-2 \text{arcsin }(a/b) \geq \frac{2\sqrt2}{\sqrt{a+1}}$ for
$b=a+1$, then we will have that inequality for all $b\geq a+1$.

Rewrite  (\ref{arcsininequality}) in terms of
$a$ and
$b=a+1$,
$$
\kappa(A)\geq \pi-2\arcsin\left(1-\frac{1}{a+1}\right)\;.
$$
Now let $t=\sqrt{\frac{1}{a+1}}$ and check (analytically or
graphically) that 
$$ 
\pi-2\arcsin(1-t^2) >(2 \sqrt{2})\;t\;.
$$

\end{proof}

In the next lemma, we call the curve $Y$ instead of $A$,  to help clarify
how the lemmas will be used later:  We  will prove Lemma
\ref{illumination} by applying Lemmas \ref{packing} and
\ref{oscillation} to subarcs $A$ of $Y$.

Suppose
$Y$ is a smooth curve in
$\R^3$, and
$x_0$ is a point some finite distance from $Y$.   The integral
\begin{equation}
\int_{y \in Y} \frac{1}{|y-x_0|^2} 
\end{equation}  can be thought of as measuring the ``illumination" of $x_0$
by $Y$.  
\begin{lemma}[Illumination and curvature]\label{illumination} Suppose $Y$ is
a smooth curve in $\R^3$, and $x_0$ is a point such that
$\forall y \in Y$, $|y-x_0|\geq 2$. Then the illumination of $x_0$ by $Y$ is
bounded by the total curvature of $Y$. More precisely,
 
\begin{equation} \label{illuminationbound}
\int_{y \in Y} \frac{1}{|y-x_0|^2} \leq c_1+c_2 \;\kappa(Y)\;,
\end{equation} where $c_1,c_2$ are universal constants independent of $Y$  {\rm (values $c_1=16$ and $c_2=43$ are sufficient).}
\end{lemma} Lemma \ref{illumination} is perhaps the most intricate part of
the paper.  Before proving it, we present four special cases.  The general
argument  does not reduce to these special cases -- rather we include them
to give an intuitive sense of why the proposition might be true (the first
four), and some of the issues one needs to confront in building a proof (the
fifth).

\subsection{Special cases for Lemma \ref{illumination}}
\subsubsection{A spiral to show the lemma is sharp in the power of $\kappa(Y)$.}\label{sharp} Let $Y$ be
the polar coordinates curve
$r=3-1/\theta$,
$\theta = 1
\ldots\Theta$.  As $\Theta$ increases, the illumnination (of $x_0=$ the
origin) is aymptotic to $\frac{1}{3}\kappa(Y)$.

\subsubsection {$Y$ is a ray} \label{raycase} Suppose $Y$ is a straight
line, starting at a point $2$ units from $x_0$ and aiming radially away from
$x_0$.  Then the line integral is just
$\int_2^{\infty}{1/s^2} \;ds = 1/2$.  
\hfill \qed 

\subsubsection {$Y$ is a straight line} Suppose $Y$ is a straight line,
infinite in both directions, and tangent to the sphere of radius $2$
centered at
$x_0$. Then 
$$\int_{y \in Y} \frac{1}{|y-x_0|^2} =
\int_{-\infty}^{\infty}\frac{1}{4+s^2}\;ds\;\;=\;\;\frac{\pi}{2}\;.
$$ If the line is a finite segment, or the minimum distance from $Y$ to
$x_0$ is $>2$, then the integral is $< \pi/2$.  \hfill \qed

\subsubsection{$Y$ is a certain kind of polygon}\label{polygoncase} Suppose
$Y'$ is a polygonal path (or closed curve)  consisting of
$e$ edges (of possibly varying lengths), such that each pair of consecutive
edges meets at a right angle.  Form a smooth curve $Y$ by replacing the
corners of $Y'$ with small quarter-circles.  Then, by the second  special
case, each edge of
$Y$ contributes  $< \pi/2$ to the illumination integral, so 
\begin{align}\nonumber \int_{y \in Y} \frac{1}{|y-x_0|^2} < e
\;\frac{\pi}{2} 
 &= \kappa(Y) \;\;\;\text{ if the polygon has endpoints, or }
\\ \nonumber &{\;\;\;\;\;}\kappa(Y) + \pi/2 \;\;\;\text{ if the polygon is
closed.}
\end{align} 
\hfill \qed

\subsubsection{$Y$ is a monotone arc}\label{monotonecase}  Suppose $Y$ is a
smooth curve, starting at a point
$y_0$ with $|y_0-x_0|=2$,  with  the property that the distance function
$|y-x_0|$ is monotone increasing on $Y$. 

For $n=2, 3, \dots$, let $B[n]$ denote the round ball of radius $n$ centered
at $x_0$, and let $S[n, n+1]$ denote the spherical shell with radii $n$ and
$n+1$. By our assumption of monotonicity, each intersection
$ Y \cap S[n,n+1]$ is a connected arc, which we denote $Y_n$. Then
$$\int_{y\in Y}\frac{1}{|y-x_0|^2}  = \sum_{n=2}^{\infty}
\int_{y\in Y_n}\frac{1}{|y-x_0|^2}\leq
\sum_{n=2}^{\infty}\frac{\ell(Y \cap S[n,n+1])}{n^2}\;.$$
 We would like to bound each of the numbers $\ell(Y \cap S[n,n+1])$, in
terms of total curvature of $Y$, somehow using   Lemma
\ref{packing}.  That lemma gives upper bounds for the lengths
$\ell(Y\cap B[n])$ in terms of total curvature, but doesn't explicitly bound
the amounts in given shells. We get around this problem by bounding (not the
illumination integral from $Y$ itself, but rather) the illumination integral
for a hypothetical curve $Y^*$ that is packed around
$x_0$ in such a way as to make the illumination integral as large as
possible subject to the constraints imposed by Lemma \ref{packing}.  (In
this intuitive discussion of the monotone case, we will continue with the
image of a ``hypothetical curve".  In the actual proof of Lemma
\ref{illumination}, we will be more rigorous.)

For brevity, let $\kappa$ denote $\kappa(Y)$.  Since $Y \cap {\rm int}\,B[2] = \emptyset$, we start with 
$Y_2 = Y\cap B[3]$. By inequality (\ref{packinginequality1}), 
$$ \ell(Y \cap B[3]) \leq 3 (\kappa + 2)\;.\\
 $$  Similarly,
\begin{align} \nonumber
 \ell(Y \cap B[4]) &\leq 4 (\kappa + 2)\;, \\ \nonumber
\ell(Y \cap B[5]) &\leq 5 (\kappa + 2)\;, \\ \nonumber &\text{etc}.
\end{align} 
  If the curve $Y$ does not actually achieve these bounds, then add extra
length (the ``hypothetical" curve $Y^*$) in each of the shells as needed to
actually reach these bounds.  Since we are adding length to the $Y$ that
already exists, the illumination integral can only increase.  Thus the
illumination for $Y^*$ is an upper bound for the illumination for $Y$.

We have
\begin{align}\nonumber
\ell(Y^* \cap B[3]) &= 3 (\kappa + 2)\;,\\ \nonumber
\ell(Y^* \cap B[4]) &= 4 (\kappa + 2)\;,\\ \nonumber
\ell(Y^* \cap B[5]) &= 5 (\kappa + 2)\;,\\ \nonumber \text{etc.}
\end{align}

Thus 
\begin{align}\nonumber
\ell(Y^* \cap S[2,3]) &= 3 (\kappa + 2)\;,\\ \nonumber
\ell(Y^* \cap S[3,4]) &= 4 (\kappa + 2)\;-\;3 (\kappa +
2)\;=\;(\kappa + 2)\;,\\
\nonumber
\ell(Y^* \cap S[4,5]) &= 5 (\kappa + 2)\;-\;4 (\kappa +
2)\;=\;(\kappa + 2)\;,\\ \nonumber \text{etc.}
\end{align}

And so,
$$  \int_{y \in Y} \frac{1}{|y-x_0|^2} \leq \int_{y \in Y^*}
\frac{1}{|y-x_0|^2} <  \frac{3 (\kappa + 2)}{2^2}
\;+\;\sum_{n=3}^{\infty}\frac{ (\kappa + 2)}{n^2} < 2\kappa+3\;.
$$ \qed

\subsection{Proof of Lemma \ref{illumination}}

We begin as we did in section \ref{monotonecase}.  For $n=2, 3, \dots$, let
$B[n]$ denote the round ball of radius $n$ centered at $x_0$, and let $S[n,
n+1]$ denote the spherical shell with radii $n$ and $n+1$. We need to bound
the total arclength of $Y$ contained in each shell $S[n,n+1]$, but we cannot
do this directly since we are not assuming monotonicity as in section
\ref{monotonecase}.  The shell-intersections might consist of long arcs, or
might consist of unions of many short arcs, as
$Y$ meanders in space, close to, or far from, 
$x_0$. 

When $Y$ is   contributing to the integral by having long arcs close to
$x_0$, we can infer curvature from Lemma
\ref{packing}. If $Y$ is contributing to the integral by oscillating in and
out from $x_0$,  we can use Lemma \ref{oscillation} to infer curvature.

To implement this plan, and handle the problem of very small oscillations,
we are going to translate the problem into discrete combinatorics.

\subsubsection{Assume finite length} If $Y$ has infinite length, express $Y$
as an increasing  union of curves of finite length.  Since the constants
$c_1, c_2$ do not depend on the curve $Y$,  we can apply inequality
(\ref{illuminationbound}) to each of these and observe that both sides of
inequality (\ref{illuminationbound}) converge appropriately.

\subsubsection{Cut $Y$ into small pieces} Pick any integer $M > \ell(Y)$,
and cut $Y$ into $M$ consecutive arcs $Y_i$ of equal length.  Let $ \epsilon
$ denote the length of each sub-arc, and note $\epsilon<1$. 

We assign to each arc $Y_i$ a label $1,2,3,\ldots$ representing the shell
$S[n,n+1]$ that (perhaps only approximately) contains $Y_i$.  Specifically,
if
$Y_i \subset S[n,n+1]$, assign label $n$. If $Y_i$ is not entirely contained
in one shell, then it must intersect a sphere $S[n]$;  because
$\epsilon < 1$, $Y_i$ can intersect at most one sphere $S[n]$;  we assign
that label $n$ to the arc.  Note that if an arc $Y_i$  carries label
$n$, then $Y_i \subset S(n-1, n+1]$ and $Y_i \cap S[n, n+1] \neq \emptyset$. 
The set of possible labels is $\left\{2, \dots, M+1\right\}$.

\subsubsection{Discretize the problem} For each integer $n$, let $\phi(n)$
denote the total number of arcs
$Y_i$ that are labeled $n$.   Thus
\begin{equation}\label{firstsum}
\int_{y \in Y} \frac{1}{|y-x_0|^2} = \sum_i \int_{y \in Y_i} \frac{1}{|y-x_0|^2}
\;<\;\sum_{n\geq 2} \phi(n) \frac{\epsilon}{(n-1)^2}\;.
\end{equation}   We also need the auxiliary function that counts the total
number of arcs that are labeled between $2$ and $n$. Define
$$ \Phi(n) = \sum_{j=2}^{n} \phi(j)\;.
$$
 We proceed as follows:
\begin{enumerate}
\item  Abstract the arc $Y$ as the string of integers $\mathcal{L}_Y=<a_1,
a_2,
\ldots, a_M>$, where $a_i$ is the shell label of $Y_i$.  
\item Show $\mathcal{L}_Y$ is constrained in certain ways.
\item Find a bound for $\Phi(n)$, in terms of $\kappa(Y)$, using Lemma
\ref{packing} and Lemma \ref{oscillation}.
\item Note that the functions $\phi$ and $\Phi$ make sense for any finite
string $\mathcal{L}$ of integers  
\item For any finite string $\mathcal{L}$ of integers $\geq 2$, define an
``energy"
$$ E(\mathcal{L}) = \sum_{n\geq 2} \phi(n) \frac{\epsilon}{(n-1)^2}\;.
$$
\item  Construct a string $\mathcal{L}^{*}$ of integers $\in \left\{2,
\dots, M+1 \right\}$ for which we know bounds on the numbers $\phi(n)$, and
for which we know
$E(\mathcal{L}) \leq E(\mathcal{L}^{*})$.

\item Find a bound for the value $E(\mathcal{L}^{*})$, which is an upper
bound for the final sum in (\ref{firstsum}), of the form we want.
\end{enumerate}

Let $\mathcal{L}_Y$ be the the string of shell labels associated to $Y$.  In
order to bound $\Phi(n)$, we first establish certain properties of
$\mathcal{L}_Y$. 

\subsubsection{Constraints on $\mathcal{L}_Y$}
 To shorten formulas in the rest of the 
proof of Lemma \ref{illumination}, we will use $\kappa$ to denote $\kappa(Y)$.

 Since the arcs $Y_i$ are listed in their order along $Y$, each intersection
$Y_i \cap Y_{i+1}$ is nonempty.  A substring such as $\langle 34434543 \rangle$ is possible. On the other hand there
cannot be a substring such as $\langle 34435543 \rangle$, because subarcs
with labels $3$ and $5$ are contained in the disjoint half-open shells $S(2,4]$ and
$S(4,6]$.  So the first constraint is:
\begin{itemize}
\item  $\mathcal{L}_Y$ must consist of contiguous labels.
\end{itemize}

For a given value of  $n$, there cannot be too many long substrings of
$\mathcal{L}_Y$ consisting of labels $\leq n$.  A substring of
$\mathcal{L}_Y$ consisting  of $q$ symbols
$\leq n$ represents a connected arc $A\subset Y$ of length $q\epsilon$
contained in the ball $B[n+1]$.  If $q$ is such that $q\epsilon \geq
3(n+1)$, then, by Lemma \ref{packing}(\ref{packinginequality2}), 
$\kappa(A) \geq 1$.  Thus, for such $q$, we can have no more than $\kappa$ such strings.  So the second constraint is:
\begin{itemize} 
\item  For each $n$, the string $\mathcal{L}_Y$ contains at most $\kappa$
pairwise disjoint substrings  of length $\frac{3(n+1)}{\epsilon}$ consisting
of integers $\leq n$.
\end{itemize}

In section (\ref{BoundPhiN}), we will apply this formula to substrings of $\mathcal{L_Y}$ with entries up through (n+1).  We also want to phrase the bound in terms
of what $\mathcal{L}_Y$ {\em cannot} contain.  Specifically,
\begin{itemize}
\item For each $n$, the string $\mathcal{L}_Y$ cannot contain $(\kappa + 1)$
pairwise disjoint substrings  of length $\frac{3(n+2)}{\epsilon}$ consisting
of integers $\leq (n+1)$.
\end{itemize}

We obtain a third constraint, this time on ``jumps".  If, in the string
$\mathcal{L}_Y$, we observe a substring $\langle  n \ldots n+1 \ldots n
\rangle$, we cannot infer any particular contribution to total curvature that is independent of $\epsilon$.  But if we
see
$\langle n \ldots n+2 \ldots n \rangle$, then we can.  Let us call  a
substring 
$\lambda = \langle n \ldots n+2 \ldots n \rangle$ of $\mathcal{L}_Y$ a {\em
jump at level $n$}.  Two jumps are {\em non-overlapping} if they are
disjoint, or meet in at most one term $a_i$ (of necessity, then, an endpoint
of each). 

An arc $Y_i$ with label $n$ has nonempty intersection with $S[n,n+1]$; an
arc with label $n+2$ intersects $S[n+2, n+3]$.  Thus if $\lambda$ is a jump
at level $n$, then the subarc of $Y$ determined by $\lambda$ has a subarc $A$  that starts at the sphere
S[n+1] and reaches as far out as some $S[b]$, $b \geq n+2$, before heading
back to end at $S[n+1]$.  By Lemma \ref{oscillation},  such an arc
contributes more than $\frac{2}{\sqrt{n+1}}$ to total curvature.  Thus we
have our third constraint:
\begin{itemize}
\item For each $n$, the string $\mathcal{L}_Y$ cannot have $\frac{1}{2} \,
\kappa \, \sqrt{n+1}$ non-overlapping jumps of level $n$.
\end{itemize}

\subsubsection{Bound $\Phi(n)$}\label{BoundPhiN}
 We now combine the constraints on substrings and jumps in
$\mathcal{L}_Y$ to get bounds on $\Phi(n)$.  
\begin{prop} \label{Phibound} For the string $\mathcal{L}_Y$, for each
$n\geq 2$,
\begin{equation} \label{PhiboundInequality}\Phi(n) < 8 \kappa
\, \frac{n^{3/2}}{\epsilon}\;+\;6\frac{n}{\epsilon}.
\end{equation}\end{prop} 
\begin{proof} Suppose, to the contrary, that for some $n$, $\mathcal{L}_Y$
does have that many symbols $2, 3, \ldots, n$.  The bound
(\ref{PhiboundInequality}) was chosen to be a simple expression  that, for
$n \geq 2$, dominates 
$$
 (\kappa + 1) \left( \frac{3(n+2)}{\epsilon} \right) \;+\; \left(
\frac{1}{2} \;\kappa \, \sqrt{n+1}\right) \left(\frac{3(n+2)}{\epsilon}\right).
$$

Visualize the string $\mathcal{L_Y}$ so that, temporarily, only the symbols
$2,3,\ldots,n$ are visible.  Parse these into  pairwise disjoint substrings
of length $\frac{3(n+2)}{\epsilon}$.  By assumption, we have (many) more than $\kappa+1$ of these substrings.  So in the actual string $\mathcal{L}_Y$, a
number of these substrings must get broken up by inserted symbols $\geq
n+1$.  Now make all the symbols $a_i=n+1$ in $\mathcal{L}_Y $ visible as
well.  These certainly can break up substrings consisting only of symbols
$2,3,\ldots,n$, but they offer no improvement on our situation of exceeding
the total curvature bound: 
 we chose the lengths of the substrings to be large enough that even if they
were made from symbols $2,3,\ldots,n+1$, they would still each be
contributing $\geq 1$ to total curvature, so we cannot have more than
$\kappa$ of these.  Thus we must have some symbols
$\geq n+2$ in the original string $\mathcal{L}_Y $ to break up a number of
the ``offending" substrings.  How many of the substrings can be broken by
inserting symbols $a_j \geq n+2 \,$?  Each offending substring that gets
broken this way  represents at least one jump at level n.  So we must have
have fewer than  $\frac{1}{2} \, \kappa \, \sqrt{n+1}$ such interruptions. 
We are assuming $\Phi(n)$ is large enough that the number of offending substrings is
greater than the number of possible interruptions plus the maxumum number we
could tolerate to be uninterrupted. We conclude that $\Phi(n)$ cannot be
that large.
\end{proof}

\subsubsection{Construct $\mathcal{L}^*$}

We have completed Step 3 of our plan, and now proceed. The functions $\phi$
and $\Phi$ make sense for abstract finite strings $\mathcal{L}$ of integers:
$$\phi(n)= \text{number of symbols } a_i \text { of } \mathcal{L} \text {
that are}\; n\;;$$
$$ \Phi(n) = \sum_{j=2}^n \phi(j)\;.
$$ 
And we can define the ``energy"
$$
 E(\mathcal{L}) = \sum_{n \geq 2} \phi(n) \frac{\epsilon}{(n-1)^2}\;.
$$

We want to construct a string $\mathcal{L}^{*}$ whose energy we can bound,
but also whose energy is larger than $E(\mathcal{L}_Y)$.  We do this by
successive modification of $\mathcal{L}_Y$. We will denote the new strings
$\L_2, \L_3, \ldots$, and denote the corresponding functions $\phi_2,
\Phi_2, \phi_3, \Phi_3$, etc. Let us also introduce notation for the bounds
in inequality (\ref{PhiboundInequality}): 
$$
\beta(n) = 8 \kappa
\, \frac{n^{3/2}}{\epsilon}\;+\;6\frac{n}{\epsilon}\;.
$$

The strings $\L_m$ will have the following properties:
\begin{itemize}
\item $\Phi_m(n)\leq \beta(n)$ for all $n$.
\item $\Phi_m(n) = \beta(n)$ for $n=2\ldots m$.
\item $E(\L_Y) \leq E(\L_m)$.
\end{itemize}

 If we take a string $\L$ and change some symbol $a_i$ to a lower integer,
that  increases $E(\mathcal{L})$.  Also, if we introduce a new additional
symbol
$a_j$ somewhere in $\mathcal{L}$, that  increases $E(\mathcal{L})$.

We begin the construction by adjoining to $\L_Y$ enough terms $a_j = (M+1)$
to raise $\Phi(M+1)$ to equal $\beta(M+1)$.  Call this string $\L_1$.  All
the other strings $\L_m$ will be obtained by changing various terms of
$\L_1$ to lower values, thus raising energy while keeping $\Phi(M+1)$
unchanged.

We know in $\L_Y$ that $\phi(2)=\Phi(2) < 8 \kappa
\, \frac{2^{3/2}}{\epsilon}\;+\;6\frac{2}{\epsilon}$. To construct $\L_2$,
change  enough 3's in $\L_Y$ to 2 to bring the number of 2's up to
$\beta(2)$.  [Note: For simplicity, we will use the bound itself, rather
than rounding up if it is not an integer.]  Changing a 3 to a 2 has no
effect on $\Phi(n)$ for $n\geq 3$.  If there are not enough 3's (i.e. if
$\Phi(3) < \beta(2)$), we change 4's to 2's.  This increases $\Phi(3)$ up to
$\beta(2)$, still $< \beta(3)$; and values $\Phi(n)$ are unchanged for $n
\geq 4$.  If there are not enough 4's, we change 5's, etc.  Continuing in
this way, we obtain a string $\L_2$ with the properties
\begin {enumerate}
\item $\phi_2(2) = \Phi_2(2) =\beta(2)$, 
\item $\Phi_2(n) \leq \beta(n)$ for $n \geq 3$, and
\item $E(\L_2) \geq E(\L_Y)$.
\end{enumerate}

We next want to make $\phi_2(3)$ large enough that $\Phi_2(3)=\beta(3)$.  So
we again change higher labels, first change 4's to 3's, then (if necessary)
5's to 3's, etc. 

We continue inductively to construct $\L_3, \ldots, \L_{(M+1)} = \L^*$, so
that
\begin{enumerate}
\item For $n=2\ldots m$, $\Phi_m(n) = \beta(n)$.
\item For $n = (m+1)\ldots (M+1)$, $\Phi_m(n) \leq \beta(n)$.
\item $E(\L_m) \geq E(\L_Y)$.
\end{enumerate}

Because we know exactly the values $\Phi^*(n)$, we can compute the values
$\phi^*(n)$ and so bound the energy.  
$$ 
\phi^*(2) = \beta(2) = 8 \kappa
\, \frac{2^{3/2}}{\epsilon}\;+\;6\frac{2}{\epsilon}\;,
$$
 and for $n \geq 3$, 
$$ 
\phi^*(n) = \Phi^*(n)-\Phi^*(n-1) = \beta(n)-\beta(n-1)   = 8 \kappa
\, \frac{n^{3/2}- (n-1)^{3/2}}{\epsilon}\;+\;\frac{6}{\epsilon} .
$$

\subsubsection{Bound $E(\L^*)$}

We bound the energy of $\L^*$ by passing to an infinite sum, so the value of
$M$ is immaterial.  Since we know the values $\phi^*(n)$, we have
\begin{align}
 \nonumber E(\L^*) &= \sum_{n=2}^{M+1} \phi^*(n) \frac{\epsilon}{(n-1)^2}
\\ &<  \nonumber 8 \kappa 2^{3/2} + 12 + 8 \kappa \sum_{n=3}^\infty
\frac{n^{3/2}-(n-1)^{3/2}}{(n-1)^2} + 6 \sum_{n=3}^\infty \frac{1}{(n-1)^2}
\\ \nonumber &< 16 + 43\kappa.
\end{align}  
\qed

\section{Proof of Theorem \ref{MainTheorem}}

As a preliminary step, rescale the knot so the thickness radius $R(K)=1$.  
This has no effect on the total curvature or on the average crossing number,
and simplifies the ratio
$E_L(K)$ to just the length, $L$. We want to show
$$
\rm{acn} (K) \leq c\cdot L\cdot\kappa(K)\;,
$$
 where $c$ is some coefficient that works for all knots. 

The average crossing number of a knot can be expressed as an integral over
the knot
\cite{FHW},  similar to Gauss's double integral formula for the linking
number of two loops. Specifically,
$$
\rm{acn}(K) = \frac{1}{4 \pi}\int_{x \in K}
\int_{y
\in K}
\frac{|<T_x\,, T_y\,, x-y>|}{|x-y|^3} \;\;,
$$
 where $T_x, T_y$ are the unit tangents at $x,y$ and $<u,v,w>$
is the triple scalar product $(u \times v)\cdot w$ of the three vectors
$u,v,w$.  

Write the double integral
 as a sum of two terms:
$$\text{Near}(K) = \int_{x \in K} \int_{\rm{arc}(x,y) \leq \pi}
\frac{|<T_x\,, T_y\,, x-y>|}{|x-y|^3} \;\;,$$ and
$$\text{Far}(K) = \int_{x \in K} \int_{\rm{arc}(x,y) \geq \pi}
\frac{|<T_x\,, T_y\,, x-y>|}{|x-y|^3} \;\;.$$

We shall analyze these contributions separately in the next two sections,
and find bounds of the form
$$ 
\nonumber \text{Near}(K) \leq b_1\;L
$$
$$
 \text{Far}(K) \leq c_1\;L + c_2\;L\;\kappa(K)\;.
$$
 where the coefficients are independent of $K$.  In each case,
we bound the inner integral and then multiply by $L$ to bound the double
integral.

Combining Near and Far, we get a bound for any smooth curve $K$ of the form
$$\text{acn}(K) \leq  aL + bL \kappa(K).$$   But if $K$ is a {\em closed}
curve, then by Fenchel's theorem, $\kappa(K) \geq 2\pi$.  Thus letting 
$c=\frac{1}{4\pi}\left(b+\frac{a}{2\pi}\right)$, we have $$\text{acn}(K) \leq c \; L  \;\kappa(K)\;.$$
Using the values of $c_1$ and $c_2$ from Lemma (\ref{illumination}) and $b_1$ from section (\ref{BoundingNearK}), we get $c \approx 3.8$.
\qed

\subsection{Bounding $\text{Near}(K)$} \label{BoundingNearK}

We shall show that the inner integral is uniformly bounded, independent of
$K$.

 For any smooth curve  with thickness radius
$R$, it is shown in
\cite{LSDR} that the curvature at each point is at most $1/R$. So in the
present situation, we know that the curvature of $K$ is everywhere $\leq 1$.

Let $\theta\to x(\theta)$ be a unit speed parametrization of $K$.  So
$x'(\theta) = T_x$ and $|x''(\theta)| \leq 1$.  We are studying points
$y$ for which $\text{arc}(x,y) \leq \pi$, so we can take for the parameter
set the interval $[0,\pi]$, with our starting point $x=x(0)$ and
$y=y(\theta)$ for some $\theta \in [0,\pi]$.  Using the  same parameter set,
let $\theta\to p(\theta)$ be an arclength preserving parametrization of the
unit semi-circle.  Since the curvature of $K$ is everywhere bounded by the
curvature of the unit circle, Schur's theorem \cite{Ch} tells us that for
each $\theta$, 
$$
 |x(\theta) - x(0)| \geq |p(\theta) - p(0)|\;, 
$$
 That is,
$$
 |y-x| \geq \sqrt{2-2\cos \theta}\;.
$$
 Thus
\begin{eqnarray} \nonumber
\frac{|<T_x\,, T_y\,, x-y>|}{|x-y|^3} \leq 
 \frac{|<T_x\,, T_y\,,
\frac{x-y}{|x-y|}>|}{2-2\cos \theta} \\ \nonumber = 
\frac{|<T_x\,, T_y\,,
\frac{x-y}{\theta}>|}{2-2\cos \theta}\;\;\frac{\theta}{|x-y|}
\end{eqnarray} Using Schur's theorem again, we have $|x-y| \geq
|p(\theta)-p(0)| = 
\sqrt{2-2\cos \theta}$. The function $\frac{\theta}{\sqrt{2-2\cos \theta}}$
is increasing on
$[0,\pi]$, with maximum value $\pi/2$.  So 
$$
\frac{|<T_x\,, T_y\,, x-y>|}{|x-y|^3} \leq \frac{\pi}{2} \;
\frac{|<T_x\,, T_y\,,
\frac{x-y}{\theta}>|}{{2-2\cos \theta}}\;.
$$
 The vectors $T_y$ and $\frac{x-y}{\theta}$ are each
first-order (in terms of $\theta$) close to
$T_x$.  Specifically, we have for $T_y$,
$$
T_y = T_x + \int_{t=0}^\theta x''(t)\;.
$$
 Since $|x''| \leq 1$, this says we can write $T_y$ as $T_x +
V$, where
$|V| \leq \theta$.  

On the other hand, 
the fundamental
theorem of calculus, applied first to $x(\theta)$ and then again to $x'(s)$, says
$$ y = x(\theta)= x(0) + \int_{s=0}^\theta
x'(s)\;ds\;=\;x(0)+\theta x'(0) +
\int_{s=0}^\theta
\int_{u=0}^s x''(u)\;du\;ds\;,
$$
so we can write
$\frac{|x-y|}{\theta}$ as $ T_x +W$, where $|W| \leq
\frac{1}{2}\theta$.

We now have
$T_x \times T_y = T_x \times \left(T_x + V\right) = T_x \times V$, which is a vector perpendicular to $T_x$ with
length $\leq \theta$.  When we take the dot product of this vector with
$T_x+W$, we just get the dot product with
$W$, so a number whose size is at most $\frac{1}{2}\theta^2$.

We now have 
$$
\frac{|<T_x\,, T_y\,, x-y>|}{|x-y|^3} \leq \frac{\pi}{4}\;\; 
\frac{\theta^2}{{2-2\cos \theta}}\;\leq
\frac{\pi}{4}\left(\frac{\pi}{2}\right)^2\;,
$$
 so the inner integral is bounded by $b_1=(2\pi)\cdot
(\frac{\pi}{4})\cdot(\frac{\pi}{2})^2$, since the points $y$ run from (what
we might denote as) $x-\pi$ to $x+\pi$.

Multiply this bound for the inner integral by $L$ to bound the double
integral.
\qed

\subsection{Bounding $\text{Far}(K)$}

As in the previous case, we shall bound the inner integral,
$$ \int_{\rm{arc}(x,y) \geq \pi}
\frac{|<T_x\,, T_y\,, x-y>|}{|x-y|^3} \;\;,$$ then multiply by $L$ to bound
the double integral. 

As before, we write the integrand as the triple scalar product of three unit
vectors, divided by $|x-y|^2$.  Since the numerator has magnitude at most
$1$, it suffices to bound 
$$\int_{\rm{arc}(x,y) \geq \pi}
\frac{1}{|x-y|^2} \;\;.$$

For any smooth curve with thickness radius
$R$, it is shown in
\cite{LSDR} that points
$x,y$ with $\text{arc}(x,y)\geq \pi R$ must have $|x-y|\geq 2R$.  So in our
situation, when $\text{arc}(x,y)\geq \pi$, we know $|x-y|\geq 2$.

Fix $x$ and let $Y=\left\{y \in K \;|\; \text{arc}(x,y) \geq \pi \right\}$. 
By Lemma \ref{illumination}, 
$$\int_Y \frac{1}{|y-x|^2} \leq c_1 + c_2 \kappa(Y) \leq c_1 + c_2
\kappa(X)\;.$$ Thus 
$$\text{Far}(K) \leq c_1 L + c_2 \kappa(X) L\;.$$
\qed

\section{Knot energies} The analysis of knot energies $E_N(K)$ and $E_S(K)$
in \cite{BS2, BS3} and the Mobius energy $E_O$ in \cite{RSBound} are similar
to the analysis of average crossing number here:  All involve bounding
``Near" and ``Far" contributions, and all rely on bounding $\int_y
\frac{1}{|y-x|^2}$ for the ``Far" part. We can use Lemma \ref{illumination}
to show that each of these energies is bounded by  [some constant, that
depends on the energy but not on $K$, times] $E_L(K) \kappa(K)$.

\section{Can the theorem be improved?} Our theorem says (throughout
this section, we will suppress coefficients)
$$\rm{acn}(K) \leq E_L(K) \;\kappa(K)\;.$$ Is it possible to lower the
exponent (from 1) on one or both of
 $E_L(K)$, $\kappa(K)$?  As noted in Example \ref{sharp}, Lemma
\ref{illumination} is sharp.  However, if we include thickness, and postulate
that the knot is long and  distributed homogeneously in space, then we can
argue heuristically that there would be a lower-order bound. This leads to
the conjecture that in fact $\rm{acn}(K) \leq E_L(K) \;\kappa(K)^{1/2}\;.$

Scale the knot $K$ so it has thickness radius $r(K)=1$. Then $E_L(K)$
is just the arclength, $L$, of $K$.  Suppose
$K$ is distributed in space so that relative to each point
$x_0
\in K$, each spherical shell $S[n,n+1]$ about $x_0$ contains on the order of 
$n^\beta$ arclength of $K$.  Here $\beta$ is constant, independent of the
choice of
$x_0$, and is a measure of the density of packing of $K$.  Fix some $x_0\in
K$. The shells run from $n=0$ to whatever value $N$ (for that $x_0$) is
needed to engulf all of
$K$.

The amount of arclength of $K$ in each shell has to include at least enough
to reach from one sphere to the other, a constant,  so $\beta\geq 0$.  On the
other hand, since $r(K)=1$, an arc (or union of arcs) of $K$ of some total length $\ell$ carries with it a proportional amount of excluded volume $=\pi \ell$. Since the volume of a spherical shell is approximately proportional to the area of a boundary sphere, we must have $\beta \leq 2$.

The total arclength $L$ of $K$ is the sum of the amounts in the shells, so
if the amount in each shell is on the order of $n^\beta$, then $L\approx N^{\beta+1}$. 

Assuming $K$ is long enough, relative to $N$, to apply Lemma \ref{packing},
we have $\kappa(K) \geq L/N \approx N^{\beta + 1}/N = N^\beta$. 

We proceed as in the proof of the main theorem: 
$\text{Near}(K)$ is bounded by a constant and we will bound  $\text{Far}(K)$ by bounding the inner Illumination integral and multiplying by $L$.  

If the amount of length of $K$ in each shell $S[n,n+1]$ is on the order of $n^{\beta}$, then 

$$\text{Illumination} \leq \sum_{n=2}^N \frac{n^\beta}{n^2} 
\approx \begin{aligned}
&\text{constant} - N^{\beta-1}\;\;\; &(0 \leq \beta < 1) \\ &\log N \;\;\; &(\beta=1) 
\\  &N^{\beta-1}\;\;\; &(1 < \beta \leq 2) \;.
\end{aligned}
$$

In the first situation $(0 \leq \beta < 1)$, we get a bound on average crossing number proportional to $L$, that is $\text{acn}(K)\leq a_1 + a_2 E_L(K)$.  This situation includes long, thin knots, as well as knots such as iterated composites of congruent curves where crossing number and ropelength grow at the same rates.
In the second situation $(\beta=1)$, we have  $\text{acn}(K) \leq a_1+a_2 E_L(K) \log \kappa(K)$.
In the third situtation $(1 < \beta \leq 2)$, we have 
$$\text{acn}(K) \leq a_1 + a_2 E_L(K) \kappa(K)^{\frac{\beta-1}{\beta}}\;.$$  
When $\beta=2$, we have the densest possible spatial packing of $K$ as in the examples \cite{B2, CKS}, where the growth rate $\text{acn}(K)\approx E_L(K)\kappa(K)^{1/2}$ is attained.

\newpage

\bibliographystyle{siam}

\bibliography{JonsEricBib_TotCurv013103}

\begin{thebibliography}{10}

\bibitem{B2}
{\sc G.~Buck}, {\em Four-thirds power law for knots and links}, Nature, 392
  (1998), pp.~238--239.

\bibitem{BO1}
{\sc G.~Buck and J.~Orloff}, {\em A simple energy function for knots}, Topology
  Appl., 61 (1995), pp.~205--214.

\bibitem{BS2}
{\sc G.~Buck and J.~Simon}, {\em Energy and length of knots}, in Lectures at
  KNOTS '96 (Tokyo), World Sci. Publishing, River Edge, NJ, 1997, pp.~219--234.

\bibitem{BS3}
\leavevmode\vrule height 2pt depth -1.6pt width 23pt, {\em Thickness and
  crossing number of knots}, Topology Appl., 91 (1999), pp.~245--257.

\bibitem{CKS}
{\sc J.~Cantarella, R.~Kusner, and J.~Sullivan}, {\em Tight knot values deviate
  from linear relations}, Nature, 392 (1998), pp.~237--238.

\bibitem{chakerian64}
{\sc G.~D. Chakerian}, {\em On some geometric inequalities}, Proc. A.M.S., 15
  (1964), pp.~886--888.

\bibitem{Ch}
{\sc S.~S. Chern}, {\em Curves and surfaces in {E}uclidean space}, in Studies
  in Global Geometry and Analysis, Math. Assoc. Amer. (distributed by
  Prentice-Hall, Englewood Cliffs, N.J.), 1967, pp.~16--56.

\bibitem{devrthickness}
{\sc Y.~Diao, C.~Ernst, and E.~J. Janse~van Rensburg}, {\em Knot energies by
  ropes}, J. Knot Theory Ramifications, 6 (1997), pp.~799--807.

\bibitem{DEJ2}
\leavevmode\vrule height 2pt depth -1.6pt width 23pt, {\em Properties of knot
  energies}, in Topology and geometry in polymer science (Minneapolis, MN,
  1996), Springer, New York, 1998, pp.~37--47.

\bibitem{devrthickness2}
\leavevmode\vrule height 2pt depth -1.6pt width 23pt, {\em Thicknesses of
  knots}, Math. Proc. Cambridge Philos. Soc., 126 (1999), pp.~293--310.

\bibitem{Fary}
{\sc I.~Fary}, {\em Sur la courbure totale d'une courbe gauche faisant un
  noeud}, Bull. Soc. Math. France, 77 (1949), pp.~128--138.

\bibitem{Fen51}
{\sc W.~Fenchel}, {\em On the differential geometry of closed space curves},
  Bull. Amer. Math. Soc. (2), 57 (1951), pp.~44--54.

\bibitem{Fox50}
{\sc R.~H. Fox}, {\em On the total curvature of some tame knots}, Ann. of Math.
  (2), 52 (1950), pp.~258--260.

\bibitem{FHW}
{\sc M.~H. Freedman, Z.-X. He, and Z.~Wang}, {\em M\"obius energy of knots and
  unknots}, Ann. of Math. (2), 139 (1994), pp.~1--50.

\bibitem{GonzalezMaddocks99}
{\sc O.~Gonzalez and J.~H. Maddocks}, {\em Global curvature, thickness, and the
  ideal shapes of knots}, Proc. Natl. Acad. Sci. USA, 96 (1999), pp.~4769--4773
  (electronic).

\bibitem{KS2}
{\sc R.~B. Kusner and J.~M. Sullivan}, {\em On distortion and thickness of
  knots}, in Topology and geometry in polymer science (Minneapolis, MN, 1996),
  Springer, New York, 1998, pp.~67--78.

\bibitem{LSDR}
{\sc R.~A. Litherland, J.~Simon, O.~Durumeric, and E.~Rawdon}, {\em Thickness
  of knots}, Topology Appl., 91 (1999), pp.~233--244.

\bibitem{Milnor53}
{\sc J.~Milnor}, {\em On total curvatures of closed space curves}, Math.
  Scand., 1 (1953), pp.~289--296.

\bibitem{MilnorTotalCurvature}
{\sc J.~W. Milnor}, {\em On the total curvature of knots}, Ann. of Math. (2),
  52 (1950), pp.~248--257.

\bibitem{Mo}
{\sc H.~K. Moffatt}, {\em The energy spectrum of knots and links}, Nature, 347
  (1990), pp.~367--369.

\bibitem{RSBound}
{\sc E.~Rawdon and J.~Simon}, {\em M\"obius energy of thick knots}, Topology
  Appl., 125 (2002), pp.~97--109.

\bibitem{RawdonIdealKnots98}
{\sc E.~J. Rawdon}, {\em Approximating the thickness of a knot}, in Ideal
  knots, World Sci. Publishing, River Edge, NJ, 1998, pp.~143--150.

\bibitem{Rawdon2000}
\leavevmode\vrule height 2pt depth -1.6pt width 23pt, {\em Approximating smooth
  thickness}, J. Knot Theory Ramifications, 9 (2000), pp.~113--145.

\bibitem{Schubert54}
{\sc H.~Schubert}, {\em \"{U}ber eine numerische {K}noteninvariante}, Math. Z.,
  61 (1954), pp.~245--288.

\end{thebibliography}

\end{document}